# Dynamic Containerized Modular Capacity Planning and Resource Allocation in Hyperconnected Supply Chain Ecosystems

Xiaoyue Liu[1], Yujia Xu[1] and Benoit Montreuil
Physical Internet Center, H. Milton Stewart School of Industrial & Systems Engineering,
Georgia Institute of Technology, Atlanta, GA 30332
Corresponding author: xliu800@gatech.edu, yxu648@gatech.edu

***Abstract:*** With the growth of data-driven services and expansion of mobile application usage, traditional methods of capacity and resource planning methods may not be efficient and often fall short in meeting rapid changes in the business landscape. Motivated by modularity, containerization, and open sharing concepts from Physical Internet (PI), this paper proposes an effective approach to determine facility capacity and production schedule to meet current and future demands by dynamically allocating Mobile Production Containers (MPCs). In this work, we develop an iterative two-stage decision making model with dynamic rolling horizon approach. The first stage is capacity planning stage, where the model determines key decisions such as project selection, facility opening periods and project-facility assignment. The second stage is resource planning stage, where the MPC allocation and relocation schedule and weekly production schedule are decided. To validate the proposed model, we conduct a case study over a modular construction supply chain focusing on the southeast US region. The results demonstrate our model not only delivers a consistent production schedule with balanced workload but also enhances resource utilization, leading to cost effectiveness.

***Keywords:*** Dynamic Facility Network Capacity Planning; Dynamic Resource Allocation and Relocation with Mobile Resources; Mobile Production Containers; Two-stage Decision Making; Physical Internet.

***Physical Internet (PI) Roadmap Fitness***: Select the most relevant area(s) for your paper according to the PI roadmaps adopted in Europe and Japan: ☒ PI Nodes (Customer Interfaces, Logistic Hubs, Deployment Centers, Factories), ☐ Transportation Equipment, ☒ PI Networks, ☐ System of Logistics Networks, ☐ Vertical Supply Consolidation, ☐ Horizontal Supply Chain Alignment, ☒ Logistics/Commercial Data Platform, ☐ Access and Adoption, ☐ Governance.

***Targeted Delivery Mode-s***: ☒ Paper, ☐ Poster, ☐ Flash Video, ☒ In-Person presentation

## 1 Introduction

According to a report from ABB Motion (2024), resource scarcity is a growing issue for 91% of industrial businesses, leading to increased costs for 37 percent of businesses, as well as supply chain disruption for 27 percent and slowdowns in production capacity for 25 percent. In addition, the COVID-19 pandemic has caused huge shortages in US labor market, which exacerbated the resource scarcity problem. The study by Ferguson (2024) indicates that the pandemic drove more than 3 million adults into early retirement as of October 2021. All of these facts highlight the urgent need to improve resource efficiency across industry, which also becomes key in achieving net zero targets.

---

[1] The authors contribute equally to this paper.

To address the problems, both industry and academia are exploring solutions to shift the conventional supply chain to a more efficient and sustainable one. In recent years, the sharing economy concept is explored to improve efficiency for modern enterprises by capitalizing on existing resources. Progress has been achieved through the implementation of this approach; however, it represents merely the initial phase, with additional research required to delve into deeper aspects. Under this global logistics challenge, Montreuil (2011) proposed the concept of Physical Internet (PI), aiming to provide an innovative way to facilitate the supply chain by utilizing modularity, containerization, and hyperconnectivity. One of the core ideas of PI is to employ the open cooperative sharing economy, emphasizing that PI infrastructure is accessible to multiple parties rather than exclusively dedicated to a single company or group of companies.

Inspired by modularity, containerization and sharing concepts from PI, we proposed an effective approach to facilitate the resource sharing of the entire supply chain ecosystem by using Mobile Production Containers (MPCs) (illustrated in Figure 1). MPC, a groundbreaking paradigm shift in containerization technology, can be considered as shipping containers with required equipment for facilities in the network, which serve as movable resources with flexibility, scalability, adaptability, and sustainability. By strategically relocating within the facility network and seamlessly integrating with potential production timelines, MPCs optimize resource utilization, facilitate agile response to evolving demands, enhance factory productivity cost-effectively, and dynamically mitigate production workload fluctuations. Additionally, their modular design and intelligent systems empower businesses to build and leverage hyperconnected supply chain ecosystems.

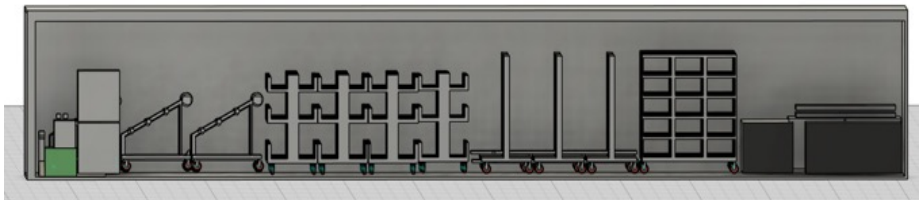

*Figure 1: Schematic Diagram for Mobile Production Containers (MPCs)*

The potential application scenarios of MPC are extensive. First, MPCs provide an agile solution for facilities facing sudden surges in demand. By swiftly relocate to these facilities, MPCs can seamlessly integrate into production lines, ensuring timely order fulfillment and optimizing resource usage. Additionally, they offer the flexibility to dynamically scale production capacity without requiring significant infrastructure investments. Second, in a globalized manufacturing scenario, MPCs offer multinational companies the ability to quickly adapt production to serve consumers worldwide. By easily deploying MPCs to small manufacturing hubs globally, transportation costs are minimized, carbon emissions decreased, and inventory management streamlined, enhancing overall efficiency. Third, in areas prone to natural disasters like hurricanes, investing in MPCs offers resilience against disruptions. By strategically deploying MPCs to nearby locations unaffected by the storm's path, production downtime is minimized. Moreover, the modular design of MPCs facilitates rapid relocation and setup, ensuring a swift response to disruptions while smoothly reintegrating original factory production processes once the storm passes.

An illustration of the relationship between capacity planning and resource deployment can be seen in Figure 2. With the input information including demand, supply, profit and cost, strategic capacity planning decisions set the maximum potential capacity of the system, while the effective and implemented capacity of the system necessary for valid production planning is determined by tactical decisions with resource deployments. By introducing the MPCs in the networks, capacity can be easily controlled dynamically by MPC allocations and relocations in facilities, taking production scheduling into consideration. This motivates the scope and focus of this research.

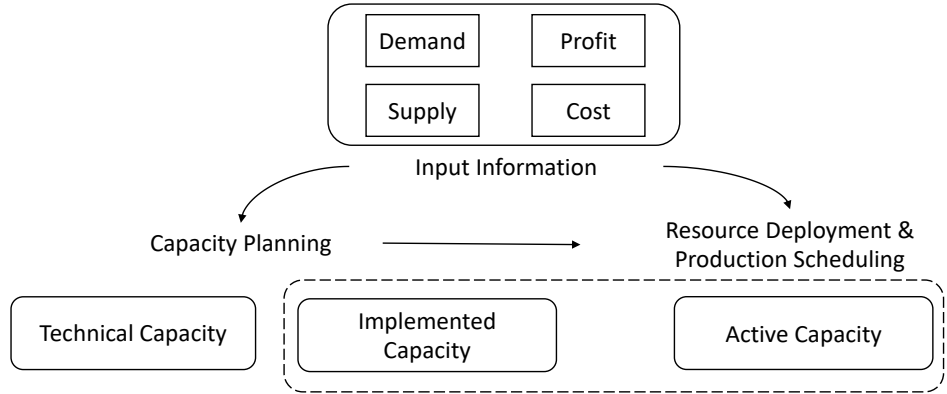

*Figure 2: Capacity planning and resource deployment*

## 2  Literature Review

A few work has been done in developing capacity planning models in the field of PI. Faugère et al. (2018, 2022) considered a hyperconnected parcel logistics system where modular smart lockers can be dynamically deployed to PI access hubs for first-mile pickup and last-mile delivery. By using those dynamic access hubs as temporal storage locations, it enables to decrease the total capacity requirements and improve the resource utilization. Moreover, Oger et al. (2021) proposed a decision support system to facilitate the strategic capacity planning more dynamic and adaptable for hyperconnected and uncertain environment. On the other hand, Liu et al. (2023, 2024) studied a dynamic modular capacity planning of logistics hub in hyperconnected relay-based transportation under the uncertainty of demand and geographical disruptions. Inspired by key concepts of modularity and uncertainty from previous works in PI, this paper considers the strategic capacity planning in hyperconnected manufacturing and production industries. By smartly deploying mobile MPCs dynamically over the production network, we can control the production throughput rate and mitigate risks sourced from fluctuate demand and unpredictable disruptions.

Dynamic and customer-oriented resource allocation is steadily expending in the literature. The Physical Internet initiative, introduced in Montreuil (2011), aims to enhance global logistics efficiency and sustainability, provides an opportunity for resource allocations and relocations in hyperconnected networks. For example, Faugère et al. (2022) works on the operations of a two-echelon synchronization problem influenced by both modular capacity reallocation and a capacity pooling recourse mechanism in urban parcel networks. Also, Xu et al. (2022) propose a dynamic workforce deployment system in hyperconnected parcel logistic hubs to match predicted demand with shifts in real-time to respond quickly to demand changes and provide guidance for centralized workforce assignment. A few papers work on the integration of capacity planning and production scheduling at a tactical level. Yao et al. (2022) presents a capacity-based tactical planning model incorporating production decisions for a multiplant multiproduct manufacturing system.

The rolling horizon (RH) approach is based on the idea that decomposing the planning horizon into several shorter, contiguous sub-periods is necessary as solving our problem over the entire horizon at once would be computationally intractable. The initial period or the first few periods for multiperiod problems are typically considered critically important since the forecasts further into the future is less reliable and more expensive, the computational burden also increase with extended planning horizons (S. Chand et al. (2002)). The global optimization model can be relaxed and broken down into multiple local optimal problems and the partial solutions in each sub-problems can be generated and assembled by utilizing this approach. Periodic RH scheduling strategy is illustrated in Figure 3(a), $n$ rolling horizon (RH) windows, each with a planning length of $T$, are employed for a problem with a planning horizon of $nT$. The choice of the time window ($T$) impacts both computation efficiency and global performance. Smaller $T$ reduces computation but sacrifices global information, potentially lowering performance (Yuan, P. et al (2018)). As illustrated in Figure 3(b), a planning horizon of $(n-1)L + T$ is decomposed into $n$ RH windows with planning length of $T$ and the rescheduling period $L$. The main difference between the two strategies is that there are overlaps between RH windows and decisions made within a RH window can be selected or be released for rescheduling in the following horizons with a certain frequency. The choice of how often to reschedule operations is important and depends on the arrival predictions and the computational costs.

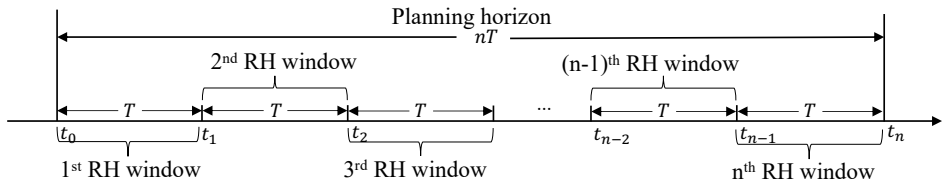

(a). Periodic RH scheduling.

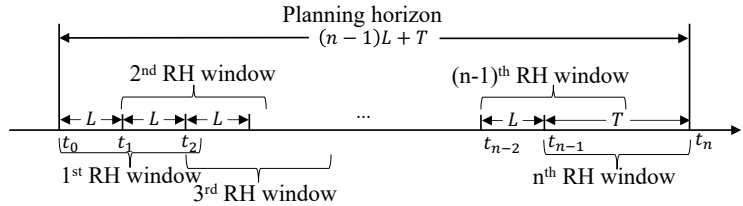

(b). Periodic RH rescheduling with selection-and-releasement.

*Figure 3: An illustration of rolling horizon (RH) strategies.*

## 3  Problem Description

In this paper, we focus on facility network capacity planning and resource allocation using MPCs. The goal of this work is to maximize the overall profit of the supply chain ecosystem by utilizing facility resources efficiently and building consistent weekly production schedules under stochastic demands. During the capacity planning stage, monthly strategic decisions involve selecting projects for each planning period, identifying new facility locations and their opening schedules, allocating accepted projects to specific facilities, and estimating the required number of MPCs to accommodate the progressive business growth. During the resource allocation stage, weekly tactical decisions involve establishing the MPC management schedule, which includes renting new MPCs, returning or relocating underutilized MPCs based on changing demand, and finalizing the production plan with smooth throughput rate.

Figure 4 illustrates the relationship between capacity planning stage and resource allocation stage. The capacity planning stage provides the information, including facility locations and status (e.g., open, close, reopen), facility-project assignment, and MPC number (i.e., capacities) at each facility, to resource allocation stage. Then, based on those provided information, the resource allocation stage will decide the weekly production schedule at each facility and update the MPC number and locations by considering relocation of MPC between multiple facilities. This approach can further reduce the MPC needed and improve the resource utilization. Finally, the unfinished demands at each facility locations and current MPC allocation are transferred back to the capacity planning stage, forming a complete loop structure.

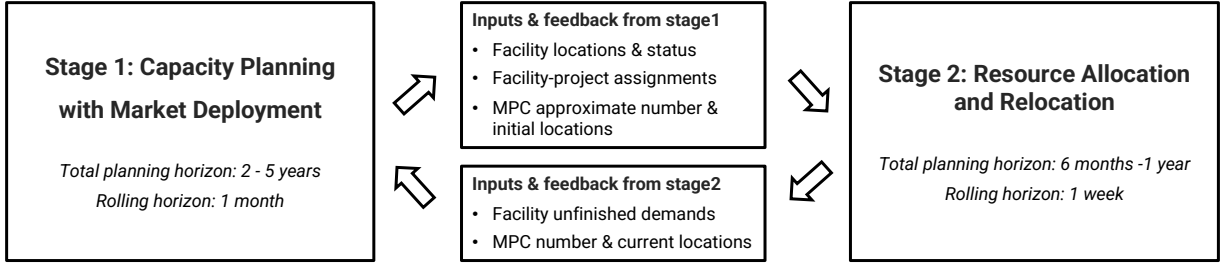

Figure 4: Relationship between Capacity Planning Stage and Resource Allocation Stage

The capacity and resource planning models are designed to operate seamlessly over time illustrated in Figure 5. At the start of each month, we conduct capacity planning with a planning horizon of 2 to 5 years. This enables us to leverage forecast demand to make long-term decisions such as negotiating project contracts, commissioning leased facilities, and preparing required resources. At the beginning of each week, we perform resource planning with a shorter planning horizon of 6 months to 1 year, which is based on the current facility network. The aim of this stage is to allocate and relocate resources, manage inventory storage, and produce weekly production. This iterative approach enables us to capture the evolving landscape of operations, facilitating agile adjustments and optimizations in response to emerging factors, including production disruptions and the emergence of new projects.

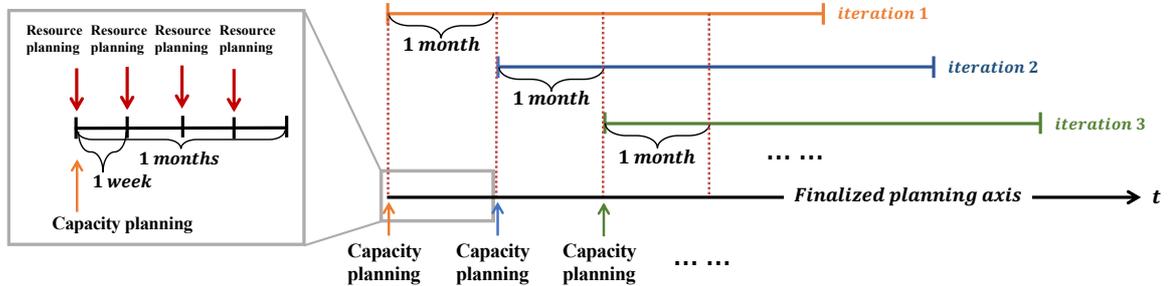

Figure 5: Timeline of the Dynamic Capacity Planning and Resource Allocation Model

## 4 Methodology

### 4.1 Capacity Planning

In this paper, we establish the following assumptions for capacity planning model. We remark that all assumptions are constructed based on realistic situations observed in a PI applied research project conducted with a real company.

- One facility can serve multiple projects.
- Each facility should be located within x miles of assigned erection sites (e.g., x=100). If a distance between facility and its assigned erection site is larger than x miles, then there will be a penalty cost based on the extra distance. The magnitude of this penalty serves as an input parameter. By tuning this parameter in line with the realistic factors of the studied case, the model can produce suitable strategies.

- Each facility has a 3-month commitment time. In other word, once a facility is opened, it should remain open for 3 months.
- Facilities can be closed after 3 months of opening based on the changing demand. In addition, each facility can be closed and reopened multiple times through the total planning horizon (i.e., 2-5 years).
- Facilities can expand or reduce its production rate by changing the number of MPCs. Specifically, each facility can lease new MPCs, return or relocate unuse MPCs during each planning period (i.e., 1 month).
- Each facility can have at most k fractals (maximum capacity limit) (e.g., k=8).
- Expand the facility capacity by one unit means leasing an additional fractal. Specifically, one additional fractal increases facility's maximal throughput rate by k (e.g., k = 2 ).
- One project can be assigned to multiple fractal facilities.
- Full truckload from facilities to erection sites and we do not consider detailed routing decisions in this study.
- The number of accepted projects, new facilities, and open facilities align with managers' strategies during each planning period (all of those can be controlled as input parameters).

Based on the provided assumptions, we develop a mixed-integer programming (MIP) model for the capacity planning stage. This model has four key decision areas: project selection, facility network planning (specifically, new facility acquisition and facility opening periods), project-facility assignment, and MPC leasing schedule. All these decisions are time dependent. The objective function aims to maximize profits within the hyperconnected supply chain ecosystem. The profits are computed as the total revenue from the accepted projects minus various costs, including new facility commissioning costs, open facility rental costs, resource commissioning and decommissioning costs, production costs, and transportation costs.

## 4.2 Resource Deployments

The dynamic containerized modular deployments in hyperconnected supply chain ecosystem are also modeled in this study. The key decisions involved include the allocation of MPCs to project activities, the scheduling of those activities, and the possible relocations of MPCs. The main assumptions for the problem are as follows:

- Decisions are made at the beginning of each week in the planning horizon.
- With the number of MPCs leased from capacity planning as an input parameter to the model, the addition or return of MPCs can be accomplished within a shorter timeframe with a detailed leasing plan.
- The initial locations of the MPCs are provided, which may either be their current positions or the most recent locations determined by the previous Resource Deployment model utilizing a rolling horizon approach.
- The relocation decisions are made at the beginning of each week.
- MPCs can only be relocated between facilities if the distance between them is within a specified distance limit.
- Transportation costs are calculated based on the shortest routing distance between facilities.
- The target start time, number of demand modules are given for each project.
- Constant erection rate, specified on a project basis, is provided as an input. One project can have one or multiple possible erection rates for consideration.
- According to current estimates, each MPC requires 2 trucks for shipment.
- Weekly production rate per MPC is assumed to remain constant.

The objective of this resource deployments model is to minimize total costs, encompassing module storage costs, production lateness costs, MPC transportation costs, MPC rental commissioning and decommissioning costs.

## 4.3 Rolling horizons

Faced with dynamic changes of demand, integrated capacity and resource planning is proposed to improve operational efficiency and flexibility utilizing a rolling horizon approach. Denote **CAP**$(t, T^l)$ and **RES**$(t, T^s)$ as two modules (or agents) for capacity planning and resource deployments at time $t$ and with planning horizon $T^l$ and $T^s$. Additionally, $S^l$ and $S^s$ denote the frequency for running capacity planning and resource deployments models, respectively. The rolling horizon procedure is detailed below in pseudocode format. To be specific, the inputs for the two proposed modules would be updated dynamically in this rolling horizon procedure.

```
ROLLING HORIZON PROCEDURE
/* Iterative process between capacity planning and resource deployments for dynamic changes */
t = 0;
WHILE t < END THEN
        IF t MOD S^l = 0 THEN
                /* Capacity planning */
                Run CAP(t, T^l);
                UPDATE Opened facilities and assigned projects;
                UPDATE Approximate number of MPCs;
        IF t MOD S^s = 0 THEN
                /* Resource deployment */
                Run RES(t, T^s);
                UPDATE Produced number of modules for each open facility;
                UPDATE MPC allocations and relocations;
                IF t MOD 28 = 0 THEN
                        UPDATE Unfinished or produced-in-advance modules with assigned projects;
                        UPDATE Number of leased or returned MPCs at facilities;
        t = t + 1;
RETURN Joint capacity and resource plans
```

# 5 Case study

## 5.1 Experimental settings

Inspired by a real-world scenario encountered within the construction industry, experiments are conducted to assess and analyze the performance of the proposed methodology. We aim to evaluate our models in a real-world scenario incorporating potential demand forecasts with rolling horizons. Given the real data of a large construction company in the US, we utilize their upcoming project forecast for capacity planning and resource deployments. The facility network and potential project locations are also assumed to be given from the company in this case study.

We set the overall planning horizon to be 20 months, the planning horizon for capacity planning to be 1 year and the planning horizon for resource deployments to be 3 months. Also, in the rolling horizon procedure, we set the frequently for running capacity planning and resource deployments models to be 1 month and 1 week, respectively.

## 5.2 Results

Figure 6 shows the required number of MPCs using three different approaches. The blue line is the baseline estimation without the proposed capacity and resource planning model, where the production rate remains fixed at 8 modules per day for each project, based on the value currently used by the company. The yellow line denotes the number of MPCs estimated through the capacity planning stage at the beginning of each month. Lastly, the green line is the required number of MPCs on a weekly basis determined by the resource allocation stage, where we consider MPC relocation and allow production date adjustments.

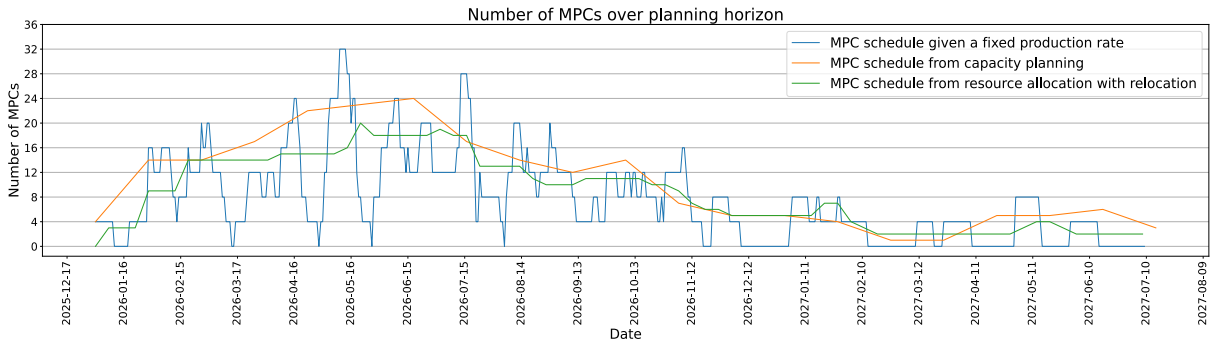

*Figure 6: Results of the required number of MPCs estimated by the fixed production rate, capacity planning model, and resource planning model with relocation*

As we can see from Figure 6, the yellow line given by the capacity planning model exhibits a smoother trend compared to the blue line from the fixed production rate scenario. This demonstrates that our model effectively balances the workload, enabling a more consistent production schedule. This is because the costs of frequently renting new fractals and returning used fractals is expensive. Since the objective of the capacity planning model is to minimize total expenses, the model will adjust the production throughput rate within each month to reduce the total number of fractals required and the frequency of changing the number of fractals. Additionally, it can be observed that the green line, representing the resource planning model with relocation, attains a lower average number of MPCs compared to the estimation derived from the capacity planning stage. This suggests that employing the relocation method can lead to further enhancements in resource utilization.

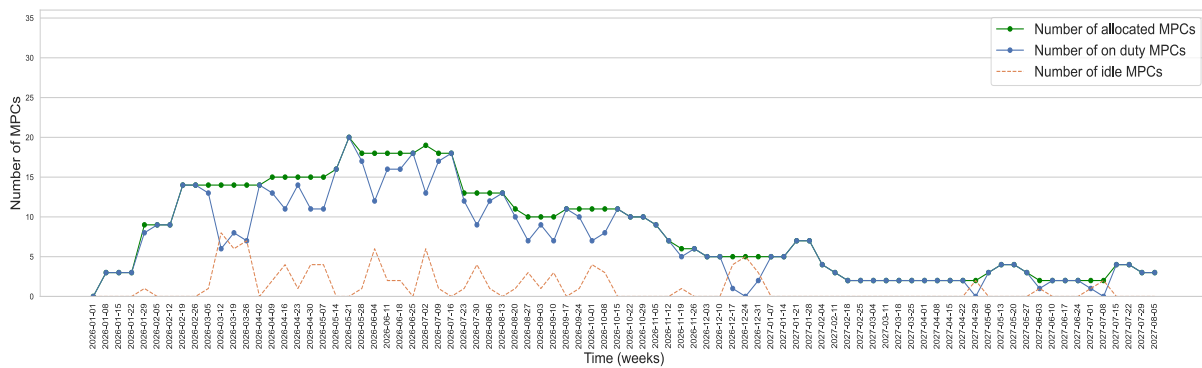

*Figure 7: Timeline depicting the number of located, idle and on-duty MPCs*

To offer a deeper insight of the MPC schedule given by the resource planning model, we present Figure 7. The green line corresponds to the same line plotted in Figure 6, which is the number of MPCs. Among them, some are on duty for module production while others are in idle, shown as blue and orange line accordingly. We can see that on duty MPC schedule from Figure 6 share the similar pattern with the blue line given a fixed production rate in Figure 7, and the number of leased MPCs has been smoothed utilizing the resource deployment model. Additionally, relocating MPCs incurs lower leasing costs for new fractals, ensuring optimal utilization of resources over time and enhancing cost efficiency within the operational framework.

# 6 Conclusion

In this paper, we present a novel rolling horizon procedure for integrated capacity planning and resource deployments in hyperconnected supply chain ecosystems. We seek to understand the value of the MPCs in capacity and resource deployments incorporating demand uncertainty and model performance under the rolling schedule procedure. By iteratively solving the proposed two models, we show that it is both computationally feasible and necessary to integrate capacity and resource deployments decisions. Dynamic changes of customer demand can be captured utilizing the integrated model with flexibility, offering management valuable insights into the optimal utilization of MPCs as demand fluctuates. For future research, stochastic programming approach considering different predictive scenarios for capacity planning under uncertainty is worth to be addressed. Additionally, distributed production planning and scheduling among multi-layer geographically separated production networks can be conducted, to determine optimal schedule of MPCs considering technical, logistical, and timing constraints.